\documentclass[review]{elsarticle}
\usepackage{lineno,hyperref}
\usepackage{amsmath}
\usepackage{amssymb}
\usepackage{theorem}
\usepackage{graphicx}
\modulolinenumbers[5]
\journal{Journal of \LaTeX\ Templates}
\theorembodyfont{\rmfamily}
\newtheorem{theorem}{Theorem}
\theorembodyfont{\rmfamily}

\theorembodyfont{\rmfamily}
\newtheorem{definition}[theorem]{Definition}
\theorembodyfont{\rmfamily}
\newtheorem{lemma}[theorem]{Lemma}
\theorembodyfont{\rmfamily}

\newtheorem{example}[theorem]{Example}

\theorembodyfont{\rmfamily}

\theorembodyfont{\rmfamily}

\newcommand{\define}{\mathrel{\hbox{$\equiv$%
 \hskip -.87em \lower .47ex\hbox{$\leftharpoondown$}}}}

\def\RR{\mathbb{R}}
\def\NN{\mathbb{N}}

\def\ds{\displaystyle}

\def\Bee{\mathcal{B}}

\def\Es{\mathcal{S}}
\def\cfun{\mathrm{1}}
\newcommand{\chFun}[1]{ {\mathrm{1}}_{#1}}
\def\vphi{\varphi}
\def\eps{\varepsilon}

\def\sfr#1#2{\frac{\ #1 \ }{\ #2 \ }}

\def\RR{{\mathbb{R}}} 
\def\NN{{\mathbb{N}} }
 
\def\ds{\displaystyle}

\def\vphi{\varphi}
\def\e{\varepsilon}

\newcommand{\intu}[1]{{\int^\uparrow_{#1}}}
\newcommand{\intd}[1]{{\int^\downarrow_{#1}}}

\def\qed{\hfill$\Box$}
\begin{document}
\begin{frontmatter}
\title{Convergence theorems for several decomposition type non-linear integrals}
\author[mymainaddress]{Ryoji Fukuda\corref{mycorrespondingauthor}}
\ead{rfukuda@oita-u.ac.jp}
\address[mymainaddress]{Oita University, 700 Dan-noharu, Oita city, 
Oita, 870-1192, JAPAN} 
\begin{abstract}
We define several types of decomposition type non-linear integrals. 
These are classified by the direction of approximation(from above or below),
the set family types (partition or covering) of simple functions, 
the coefficient signature (non-negative or signed),
and cardinal number of terms of simple functions(finite or countable infinite).
We will compare these integrals 
considering the monotone increasing/decreasing convergence theorems.
\end{abstract}
\end{frontmatter}
\section{Introduction}
This is an English translation of 
``Comparison of Decomposition Type Nonlinear Integrals Based on the Convergence Theorem''
\cite{FHO2020J}.
We define several decomposition type non-linear integrals.
The view points are the direction of approximation,
the set family types of simple functions, 
the coefficient signature,
and cardinal number of terms of simple functions.
We will give some sufficient conditions for convergence theorems:
monotone increasing convergence theorems,
monotone increasing convergence theorems,
and uniform convergence theorems.

\section{Classes of Simple functions  and  Definitions of Integrals}
We will give some concepts and notations. 
 Throughout the paper, $(X,\Bee)$ denotes a measurable space.
 $X$ is non-discrete set and $\Bee$ is a $\sigma$-algebra.
 We call a set function $\mu$ ($\Bee\to\RR^+$) ``a monotone measure'' 
if $\mu(\emptyset)=0$ and $\mu(A)\leq \mu(B)$ if $A\subset B$. 
 We assume that all monotone measures $\mu$ satisfies continuity from above and below, that is:
\[
		A_n\nearrow A, \mbox{ or } A_n\searrow A \mbox{ as } n\to\infty \quad
	\Rightarrow \quad 
		\mu(A_n) \to \mu(A)\mbox{ as } n\to\infty.
\]
Let $\vphi$ be a simple function expressed by
$\vphi(x)=\sum_k a_k \chFun{A_k}$,
where $a_k\in \RR$ and $A_k\in\Bee$ for each $k$.
The summation may be finite or infinite.
For this simple function, we define the {\it basic sum} $\mu(\vphi)$ 
of $\vphi$ with respect to $\mu$ by
\[
	\mu(\vphi)=\sum_k a_k \mu(A_k).
\]
We assume that the series converges absolutely
when the summation is infinite.
We do not assume the additivity for a monotone measure $\mu$,
then the basic sums are not the same 
among a family of simple functions 
which are same as functions.
Hence, we have to distinguish simple functions when 
sequences of pairs of a real number and a measurable set are not the same
when we consider the basic sums.
\begin{definition}\label{sFunClassInt}
We define 8 simple function families as follows.
Let $\mu$ be a monotone measure.
\begin{eqnarray*}
   \Es^{P+}&=&\{\{(a_k,A_k)\}_{k=1}^n: n\in\NN, a_k \geq 0, \ \{A_k\}_k \mbox{ is a partition of } X \} \\
   \Es^{P\pm}&=&\{\{(a_k,A_k)\}_{k=1}^n: n\in\NN, a_k\in\RR, \ \{A_k\}_k \mbox{ is a partition of } X \} \\
   \Es^{P+}_\mu&=& \{\{(a_k,A_k)\}_{k=1}^\infty:  a_k\geq 0, \sum_k a_k\mu(A_k)<\infty,\\
					&&			\hspace{12em} \{A_k\}_k \mbox{ is a partition of } X \} \\
   \Es^{P\pm}_\mu&=& \{\{(a_k,A_k)\}_{k=1}^\infty:  a_k\in\RR, \sum_k |a_k|\mu(A_k)<\infty,\\
					&&			 \hspace{12em} \{A_k\}_k \mbox{ is a partition of } X \} \\
   \Es^{C+}&=&\{\{(a_k,A_k)\}_{k=1}^n: n\in\NN, a_k\geq 0, \ \{A_k\}_k \mbox{ is a covering of } X \} \\
   \Es^{C\pm}&=&\{\{(a_k,A_k)\}_{k=1}^n: n\in\NN, a_k\in\RR, \ \{A_k\}_k \mbox{ is a covering of } X \} \\
   \Es^{C+}_\mu&=& \{\{(a_k,A_k)\}_{k=1}^\infty:  a_k\geq 0, \sum_k a_k\mu(A_k)<\infty, \\
					&&\hspace{12em}						 \{A_k\}_k \mbox{ is a covering of } X \} \\
   \Es^{C\pm}_\mu&=& \{\{(a_k,A_k)\}_{k=1}^\infty:  a_k\in\RR, \sum_k |a_k|\mu(A_k)<\infty,\\
					&&\hspace{12em}						 \{A_k\}_k \mbox{ is a covering of } X \} \\
\end{eqnarray*}
A simple function $\vphi= \{(a_k,A_k)\}$ is a sequence of pairs of a real number and a measurable set.
We always identify $\vphi$ with the function
\[
	\vphi(x) =\sum a_k \cfun_{A_k}(x).
\]
For a family of simple functions $\Es$ and a measurable function $f$ on $(X,\Bee)$,
we define the following families of simple functions.
\begin{eqnarray*}
   L(\Es,f) &=& \{ \phi \in\Es, \phi (x) \leq f(x),\ \forall x\in X\}. \\
   U(\Es,f) &=& \{ \phi \in\Es, \phi (x) \geq f(x),\ \forall x\in X\}.
\end{eqnarray*}
Then, we define decomposition type integrals as follows.
\[
	\intu{\Es} f d\mu= \sup\{\mu(\vphi) : \vphi\in L(\Es,f) \},\quad
	\intd{\Es} f d\mu= \inf\{\mu(\vphi) : \vphi\in U(\Es,f) \}.
\]
$\intu{\Es^{P+}}$ is a Pan integral(\cite{Yang1985}),
$\intu{\Es^{C+}}$ is a SD integral(\cite{MesiarPap2015}),
$\intd{\Es^{P+}}$ is a concave integral(\cite{Lehrer2009}) and
$\intd{\Es^{C+}}$ is a convex integral(\cite{MesiarPap2015}).
\end{definition}
We formally defined several integrals, however, 
the simple function classes $\Es^{c,\pm}$ and $\Es^{c,\pm}_\mu$
are too wide to define the above integrals for standard measurable functions.

Next we will explain some basic properties of these integrals.
\begin{lemma} \label{monotone}
Let $\mu$ be a monotone measure, 
$f,g$ be measurable functions on $(X,\Bee)$, 
and $c$ be a positive constant.
Then, for each $\int=\intu{\Es}, \intd{\Es}$, 
and $\Es$:one of the simple function classes defined in Definition \ref{sFunClassInt},
\begin{description}
\item[(a)] $f\leq g$ implies $\ds \int f d\mu \leq \int g d\mu,$
\item[(b)] $\ds \int c f d\mu = c \int f d\mu.$
\end{description}
\end{lemma}
Proof.\quad
(a) $f\leq g$ implies
\[
	\{ \vphi : \vphi \leq f \} \subset 	\{ \vphi : \vphi \leq g\},\quad
	\{ \vphi : \vphi \geq f \} \supset 	\{ \vphi : \vphi \geq g\}.
\] 
Then
\[
	\sup \{ \mu(\vphi) : \vphi \leq f \} \leq \sup	\{ \mu(\vphi) : \vphi \leq g\}, 
\]
\[
	\inf \{ \mu(\vphi) : \vphi \leq f \} \geq	\inf\{ \mu(\vphi) : \vphi \geq g\}.
\] 
These conclude the proof of (a).

(b) This relation can be easily obtained, from the fact that  
\begin{description}
\item[(c)] $c\vphi\in\Es$, $\mu( c \vphi)=c\mu(\vphi)$,
\item[(d)] $f\leq \ (\geq) \vphi$ if and only if $cf \leq \ (\geq) c \vphi$.
\end{description}
\mbox{\quad}\qed
\begin{lemma}\label{decIneq1}
Let  $\mu$  be a monotone measure, $\delta>0$ be a positive number, and $f$ be a nonnegative function on $(X,\Bee)$. 
When $\Es=\Es^{P+}, \Es^{P+}_\mu$,
\[
	\intu{\Es} f + \delta \cfun_X d\mu \leq 
	\intu{\Es}  f d\mu +\delta \intu{\Es} \cfun_{X} d\mu.
\]
\end{lemma}
Proof.\quad
Let $\vphi=\sum_k a_k\cfun_{A_k} $ be an element of $L(\Es,f+\delta\cfun_X)$.
We may assume that $a_k$ not less than  $\delta$, 
since $f+\delta\cfun_X$ is not less than $\delta$.
Then,
\[
	\vphi_1=\sum_k (a_k-\delta )\cfun_{A_k} \in L(\Es,f), \quad 
	\vphi_2=\sum_k  \cfun_{A_k} \in L(\Es,\cfun_X).
\]
By the definition of $\vphi_1$ and $\vphi_2$, 
\[
	\mu(\vphi) =\mu(\phi_1+\delta\vphi_2) =\mu(\phi_1)+\delta\mu(\vphi_2).  
\]
Therefore,
\begin{eqnarray*}
		\intu{\Es} (f+\delta\cfun_{X}) d\mu
&=& 	\sup \{ \mu(\vphi) ;\vphi\in L(\Es,f+\delta\cfun_X) \} \\
&=& 	\sup \{ \mu(\vphi_1)+\delta\mu(\vphi_2) ;\vphi_1\in L(\Es,f),\ \vphi_2\in L(\Es,\cfun_X) \} \\
&\leq&\sup \{ \mu(\vphi_1): \vphi_1\in L(\Es,f)\} +\delta \sup\{\mu(\vphi_2); \vphi_2\in L(\Es,\cfun_X) \} \\
&=&	\intu{\Es} f d\mu + \delta \intu{\Es} \cfun_X d\mu. 
\end{eqnarray*}
\mbox{\quad}\qed
\begin{lemma} \label{decIneq2}
Let $\mu$  be a monotone measure, $f,g$ be nonnegative functions on $(X,\Bee)$, 
and $A$ is a $\Bee$-measurable set.
	When $\Es=\Es^{P+}, \Es^{P+}_\mu$, 
\begin{description}
\item[(a)] 
\[
	\intu{\Es} f d\mu \geq \intu{\Es}  f\cfun_{A} d\mu +\intu{\Es} f\cfun_{A^c} d\mu, 
\]
\item[(b)] 
\[
	\intd{\Es} f d\mu \leq \intd{\Es}  f\cfun_{A} d\mu +\intd{\Es} f\cfun_{A^c} d\mu. 
\]
\end{description}
	When $\Es=\Es^{C+}, \Es^{C+}_\mu$, 

\begin{description}
\item[(c)] 
\[
	\intu{\Es} f+g d\mu \geq \intd{\Es} f d\mu +\intd{\Es} g d\mu, 
\]
\item[(d)] 
\[
	\intd{\Es} f+g d\mu \leq \intd{\Es} f d\mu +\intd{\Es} g d\mu. 
\]
\end{description}

\end{lemma}
Proof.\quad
(a)\quad 
We consider simple functions $\vphi_1 \in L(\Es,f\ \cfun_A)$and $\vphi_2\in L(\Es,f\ \cfun_{A^c})$,
and assume that these are expressed by
\[
 	\vphi_1 = \sum_k b_k \cfun_{B_k},\quad
 	\vphi_2 = \sum_k b'_k \cfun_{B'_k}.
\]
Then $B_k\cap A^c\not=\emptyset \Rightarrow b_k=0$, and 
$B'_k\cap A \not=\emptyset \Rightarrow b'_k=0$,
since $\cfun_A=0$ on $A^c$ and  $\cfun_{A^c}=0$ on $A$. 
After removing sets with $b_k=0$ or $b'_k=0$, the family $\{B_k\} \cup \{B'_k\}$ is a disjoint family.
Hence, $\vphi =\sum_k b_k \cfun_{B_k}+ \sum_k b'_k \cfun_{B'_k}\in \Es^{C+}, \Es^{C+}_\mu$ and
$ \vphi \in L(\Es,f )$. This implies
\begin{eqnarray*}
			\intu{\Es} f d\mu
&=& 		\sup\{ \mu(\vphi) : \vphi \in L(\Es,f)\}   \\ 
&\geq& 	\sup\{ \mu(\vphi_1) : \vphi_1 \in L(\Es,f \cfun_{A}) \} + \sup\{ \mu(\vphi_2) : \vphi_2 \in L(\Es,f \cfun_{A^c}) \}\\
&= & 	\intu{\Es} f \cfun_A d\mu+\intu{\Es} f \cfun_{A^c} d\mu.
\end{eqnarray*}
(b)
Let $\vphi_1 $and $\vphi_2$ be simple functions with $\phi_1 \in U(\Es,f \cfun_{A})$, 
$\phi_2 \in U(\Es,f \cfun_{A^c})$.
When, $B_k$ is replaced by $B_k\cap A^c$ and 
$B'_k$ is replaced by $B'_k\cap A$, the following properties still hold.
\[
	\vphi_1 \in U(\Es,f \cfun_{A}), \quad 
	\vphi_2 \in U(\Es,f \cfun_{A^c}).
\]
Then, $\vphi=\vphi_1+\vphi_2 \in U(\Es,f)$, and this implies
\begin{eqnarray*}
			\intd{\Es} f d\mu
&=& 		\inf\{ \mu(\vphi) : \vphi \in U(\Es,f)\}   \\ 
&\leq& 	\inf\{ \mu(\vphi_1) : \vphi_1 \in U(\Es,f \cfun_{A}) \} + \inf\{ \mu(\vphi_2) : \vphi_2 \in U(\Es,f \cfun_{A^c}) \}\\
&=& 	\intd{\Es} f \cfun_A d\mu+\intd{\Es} f \cfun_{A^c} d\mu.
\end{eqnarray*}
(c) (d) \quad
Let $\vphi_1,\vphi_2$ be simple functions with $\vphi_1\in L(\Es,f)\ (U(\Es,f))$ and  $\vphi_2\in L(\Es,g)\ (U(\Es,g))$.
By the definition of $\Es^{c+},\ \Es^{c+}_\mu$, $f+g\in  L(\Es,f+g)\ (U(\Es,f+g))$.
This implies that
\begin{eqnarray*}
&&		\intu{\Es} f+g d\mu,\quad \left( \intd{\Es} f+g d\mu \right)\\  
&=& 	\sup\ (\inf\ )\{ \mu(\vphi): \vphi\in L(\Es,f+g) \ (\ U(\Es,f+g)\ )\} \\
&\geq (\leq)& \sup \ (\inf\ ) \{ \mu(\vphi_1): \vphi_1\in L(\Es,f) \ (\ U(\Es,f)\ )\} \\
&&	+ \sup\ (\inf\ )\{ \mu(\vphi_2): \vphi_2\in L(\Es,g) \ (\ U(\Es,g)\ )\} \\
&=& \intu{\Es}\ \left( \intd{\Es} \right)\ f d\mu +\intu{\Es}\ \left( \intd{\Es} \right)\  g d\mu
\end{eqnarray*}
This concludes the proof.\qed
\section{Uniform convergence theorem.}
First, we consider the uniform convergence theorem for Pan integral.
\subsection{Uniform convergence theorem for Pan integral}

\begin{lemma}\label{preUcPan}
Let $\mu$ be a monotone measure
and $f$ be a measurable function. 
Assume that $\ds\int \cfun_X d\mu =M<\infty$.
\begin{description}
\item[(a)]\quad
When $\Es=\Es^{P+},\ \Es^{P+}_\mu$ and $f$ is nonnegative,
\[
	\intu{\Es}  f d\mu - \delta M 
	\leq\intu{\Es}  (f - \delta) \vee 0 d\mu 
	\leq\intu{\Es}  (f + \delta)  d\mu 
	\leq \intu{\Es}  f d\mu + \delta M 
\]
for any $\delta>0$. 
\item[(b)]\quad
When $\Es=\Es^{P\pm},\ \Es^{P\pm}_\mu$,
\[
	\intu{\Es}  f d\mu - \delta M 
	\leq\intu{\Es}  (f - \delta)  d\mu 
	\leq\intu{\Es}  (f + \delta)  d\mu 
	\leq \intu{\Es}  f d\mu + \delta M 
\]
for any $\delta>0$. 
\end{description}
\end{lemma}
Proof.\quad
	We will prove the third inequality for (a) and (b), that is, 
	$\Es=\Es^{P+},\ \Es^{P+}_\mu,\ \Es^{P\pm},\ \Es^{P\pm}_\mu$.
	We assume that $f\geq 0$ if $\Es=\Es^{P+},\ \Es^{P+}_\mu$.
	
	Fix an arbitrary $\eps>0$. Then there exists$\vphi\in\Es$ such that 
\[
		\vphi\leq f+\delta,\quad
		\mu(\vphi) \geq \intu{\Es} (f+\delta ) d\mu -\eps .
\]
Using the representation $\vphi=\sum_k a_k\cfun_{A_k}$,
we may assume that $a_k\geq \delta$ if $\Es=\Es^{P+},\ \Es^{P+}_\mu$.
We define new simple function 
\[
	\psi =\sum_k (a_k-\delta)\cfun_{A_k} \leq f,  
\]
Remark that the coefficients $(a_k-\delta)$ are non-negative when $\Es=\Es^{P+},\ \Es^{P+}_\mu$.
In any cases, we have 
\begin{eqnarray*}
	\intu{\Es} (f +\delta) d\mu -\eps 
&\leq & \mu(\vphi) \\ 
&\leq & \mu(\psi) +	\mu\left( \sum_k \delta \cfun_{A_k} \right)\\
&\leq & \int f d\mu +\delta \int\cfun_X d\mu .\\
&= 	& \int f d\mu +\delta M.
\end{eqnarray*}
This implies that
\[
		\int (f+\delta ) d\mu \leq \int f d\mu +\delta M,
\]
since $\eps$ is any positive number.  

The second inequality in (a) and (b) are obvious, and we will prove the first one.
We consider the case $\Es=\Es^{P+},\ \Es^{P+}_\mu$,
For any $\e>0$, there exist $\vphi\in\Es$ with 
\[
		\vphi\leq f,\quad
		\mu(\vphi) \geq \intu{\Es} f d\mu -\eps .
\]
Using the representation $\vphi=\sum_k a_k\cfun_{A_k}$,
define a new simple function
\[
	\psi=	\sum_k \{ (a_k-\delta)\vee 0 \} \cfun_{A_k} \leq (f-\delta)\vee0,  
\]
Hence, $\psi\in L(\Es, (f-\delta)\vee0)$.
\begin{eqnarray*}
			\intu{\Es} f d\mu -\eps 
&\leq& 	\mu(\vphi) \\  
&\leq&	\mu(\psi) + \delta\sum \mu(A_k)\\
&\leq&  	\intu{\Es} (f-\delta)\vee 0 d\mu+ \delta M.
\end{eqnarray*}
Then we have
\[
			\intu{\Es} f d\mu -\delta M \leq \intu{\Es} (f-\delta)\vee 0 d\mu.
\]
Proof of the first inequality, for the case $\Es=\Es^{P\pm},\ \Es^{P\pm}_\mu$,
is parallel with the above proof.
\ \qed

\begin{theorem}\label{ucPan}
Let $\mu$ be a monotone measure, and
$\Es=\Es^{P+},\Es^{P\pm},\Es^{P+}_\mu,\Es^{P\pm}_\mu$.
Assume that $\ds M=\intu{\Es} \cfun_X d \mu<\infty$ . 
Then, if a sequence of nonnegative measurable functions $\{ f_n \}$ 
converges to $f$ uniformly and $\intu{\Es} f d\mu<\infty$, 
moreover, we also assume that $f$ and $f_n$ ($n\in\NN$) are non-negative if $\Es=\Es^{P+},\Es^{P+}_\mu$,
then
\[
	\lim_{n\to\infty} \intu{\Es} f_n d\mu= \intu{\Es} f d\mu.
\]
\end{theorem}
Proof.\quad
For any  $\delta>0$, there exists $n_\delta\in\NN$ such that,
\[
		f(x)-\delta \leq f_n(x) \leq f(x)+\delta
\]
for any $n\geq n_\delta$.
When $\Es=\Es^{P+},\Es^{P+}_\mu$,
\[
		(f(x)-\delta)\vee 0 \leq f_n(x) \leq f(x)+\delta
\]
for any $n\geq n_\delta$. Using Lemma \ref{preUcPan}, we have
\[
		\intu{\Es} f_n d\mu 
		\in \left[\intu{\Es} f d\mu -\delta M, \intu{\Es} f d\mu +\delta M \right]
\]
for any $n\geq n_\delta$.
Hence $\intu{\Es} f_n d\mu$ converges to $\intu{\Es} f d\mu$ as $n\to\infty$.\qed

\subsection{Uniform convergence theorem for concave integral}
We consider the case $\Es = \Es^{c+},\Es^{c\pm},\Es^{c+}_\mu,\Es^{c\pm}_\mu$.
The next example illustrates that conditions for the uniformly convergence theorem 
are different for concave integral.

\begin{example}\label{exaCav}
Set $X=\NN_0=\{0,1,2,\ldots\}$, and a monotone measure $\mu$ is defined by 
\[
	\mu(A) \quad = \quad \left\{\quad
	\begin{matrix} 
		0 \ ,&\quad \mbox{$A$ is one point set or $0\not\in A$,}\\
		1 \ ,&\quad |A| >1 \mbox{ and } 0\in A.
	\end{matrix}
	\quad \right.
\]
For each $n\in\NN$ we define a function $f_n$ ($n\in \NN)$ as follows. 
\[ 
	f_n(k) = \left\{\quad \begin{matrix}
 								1             \ ,&& \quad\  k=0, \\
				 \ds \sfr{1}{n} \ ,&&\quad\ \mbox{otherwise.} 
\end{matrix}\right. \]
Then the following properties (a) $\sim$ (d) hold.
\begin{description}
\item[(a)] \quad $\mu$ is continuous from below, and is not continuous from above. 
\item[(b)] \quad When $\Es=\Es^{C+},\Es^{C+}_\mu$,
\[
	\intu{\Es}  \cfun_{\NN_0}  d\mu< \infty.
\]
\item[(c)] \quad $f_n\searrow \cfun_{\{0\}}$ uniformly.
\item[(d)] \quad For all $n\in\NN$,
\[ 
	\intu{\Es}  f_n  d\mu \cfun_{\NN_0}=1
	\not= \intu{\Es}  \cfun_{\{0\}} d\mu=0.
\]
\end{description}
\end{example}
Proofs and Comments. \quad
(a)  Let $\{ A_n \}$ be a sequence of measurable sets satisfying $A_n \nearrow A$. 
If $\mu(A)=1$, $A$ contains $0$ and other one point $a_0$.
Then $a_0, 0\in A_n$ for large enough $n\in \NN$.
This implies that  $\mu(A_n)=1$ and $\mu$ is continuous from below.

Set $A_n=\{0, n,n+1,\ldots\}$, then 
\[
	\bigcap_{n=1}^\infty A_n=\{ 0 \},\ \mu(A_n)=1,\ \mu(\{ 0 \})=0,
\]
This prove the discontinuity of $\mu$ from above.

(b) Let $\vphi $ be an element of $L(\Es,\cfun_{\NN_0})$,
$\vphi =\sum_k b_k\cfun_{B_k} \leq 1$. Then,
 \[
		\sum_{0\in B_k} b_k \leq 1,
\]
since $0\not\in B_k$ implies $\mu(B_k)=0$. 
\begin{equation}
	\mu(\vphi) =\sum_k b_k \mu(B_k)\leq \sum_{0\in B_k} b_k \leq 1.
	\label{leq1} 
\end{equation}
Hence, $\intu{\Es} \cfun_{\NN_0} d\mu <\infty$.

(c) obvious. 

(d)\quad
Set $B_k=\{0,k\}$,  
and $\vphi_n = \sum_{k=1}^n \sfr{1}{n} \cfun_{B_k}$, then, 
\[ 
		\vphi_n \leq f_n,\ \mu(\vphi_n )=1.
\]
This implies,
\[ 
	\intu{\Es} f_ n d\mu\geq 1.
\]
By the inequality (\ref{leq1})  and $f_n\leq \cfun_{\NN_0}$,  
\[ 
	\intu{\Es} f_ n d\mu\leq \intu{\Es}  \cfun_{\NN_0} d\mu\leq 1.
\]
Hence, $\intu{\Es} f_ n d\mu =1$. 

Let $\vphi \in L(\Es,\cfun_{\{0\}})$, and $\vphi= \sum a_k \cfun_{A_k}$.
$A=\emptyset$ or $A=\{ 0 \}$, then,
the summation is single $\vphi= a_1 \cfun_{\{ 0\}}$
and $a_1\leq 1$.
Therefore, we have $\intu{\Es} \cfun_{\{0\}} d\mu=0$, since $\mu(\{0\})=0$.

\[
	\intu{\Es}  \cfun_{\{0\}} d\mu=0 \not= \lim_{n\to\infty} \intu{\Es}  f_n d\mu.
\]
Thus, the uniform convergence theorem is not valid.\qed

\begin{theorem}
Let $\mu$ be a monotone measure, 
$\{ f_n \}$, $f$ be non-negative measurable functions,
and $\Es=\Es^{c+},\Es^{c+}_\mu$.
Assume that $\mu$ is continuous from below,
\[ 
	\inf_{n\in \NN,x\in X} f_n(x)=a>0,\quad
	\lim_{n\to\infty}\sup_{x\in X} |f_n(x)-f(x)|=0,
\]
and
\[
	\int f_n d\mu,\int f d\mu<\infty	
 \]
 Then, 
\[
	\int f_n d\mu \to \int f d\mu.
\]
\end{theorem}
Proof.\quad
By the assumption that $f(x)\geq a$ ($x\in X$).
 For any $\delta\in (0,1)$
there exists $\eps>0$ such that
\[
	(1-\delta) f (x)< f(x)-\eps <f(x)+\eps<(1+\delta)f(x)
\]
for any $x\in X$.
Using the uniform convergence to $f$,
there exists $N\in\NN$ such that
$|f_n(x)-f(x)|<\eps$ for any $x\in X$ and $n\geq N$.
This implies
\[
	(1-\delta) f < f_n<(1+\delta)f
\]
Therefore,
\[
	(1-\delta)\intu{\Es} f d\mu\leq \intu{\Es} f_n d\mu
	\leq (1+\delta)\intu{\Es} f d\mu.
\]
Thus, we conclude the proof by $\delta\to 0$.  \qed

\section{Monotone Convergence Theorem}
In this section, we discuss about monotone increasing and decreasing
convergence theorems, these properties are deeply 
connected with the approximation direction used in
the definition of integrals.
 
\subsection{Monotone increasing convergence theorem for $\intu{\Es}$} 
For simple function families $\Es=\Es^{p+},\Es^{p\pm},\Es^{p+}_\mu,\Es^{p\pm}_\mu,\Es^{c+},\Es^{c+}_\mu$,
we will prove the monotone increasing convergence theorem for $\intu{\Es}$
using an essentially same method.
For the classes $\Es^{c\pm},\Es^{c\pm}_\mu$,
$L(\Es,f)$ or $U(\Es,f)$ are too wide and  
the corresponding integrals do not make sense.
Then, we do not treat these integrals.

When $\Es$ is a family of infinite sum, 
we need the following properties,
which can be easily proved using the dominated convergence 
theorem (see for example  \cite{TelenceTaoBook}).

\begin{lemma}\label{sum_dct}
Let $\{ a_k\}_k,\  \{ x_k\}_k$ be real sequences, 
$\{\ \{ x^{(n)}_k \}_k\ \}_{n\in \NN}$ be a sequence of real sequences.
We assume that
\begin{description}
	\item[(a)] $a_k\geq 0$ for any $k\in\NN$ and
		$\ds \sum_k a_k<\infty $.
	\item[(b)] $|x^{(n)}_k|\leq a_k $ for any $n,k\in\NN$.
	\item[(c)] $\ds \lim_{n \to \infty} x^{(n)}_k= x_k$ for any $k\in \NN$.
\end{description}
Then,
\[
		\lim_{n\to \infty} \sum_k x^{(n)}_k  =\sum_k x_k . 
\]
\end{lemma}
\ \ \qed

In the case $\Es$ consist of non-negative functions,
we have the following theorem.
\begin{theorem}\label{incCt}	
	When $\Es=\Es^{p+},\Es^{p+}_\mu,\Es^{c+},\Es^{c+}_\mu$,\ 
	$\intu{\Es}$ satisfies monotone increasing convergence theorem.
	That is, increasing sequence $\{f_n\}$ of non-negative functions
	converges to $f$.
	Then we have
\[
	\lim_{n\to\infty} \intu{\Es} f_n d\mu= \intu{\Es} f d\mu.
\] 
\end{theorem}
Proof.\quad
Set $M=\int f d\mu$.
Let $\eps>0$ be an arbitrary positive number,
and set $M'= M-\eps$ if $M<\infty$.
If $M=\infty$, let $M'$ be any positive number. 
Then, we can select $\vphi \in L(\Es,f)$ with
\[
	 \mu(\vphi) > M'.
\]
For $\delta>0$, we define 
\[
	A^{(\delta)} _n = \{ x| f_n(x)\geq f(x)(1-\delta) \}.
\]
Then, $A^{(\delta)}_n \nearrow X$ as $n\to \infty$.
We define 
\[
	\vphi_n = (1-\delta)\vphi \cfun_{A^{(\delta)}_n}.
\]
Then, $\vphi_n\in\Es$ and $\vphi_n \leq f_n$.  
By Lemma \ref{sum_dct}, we have
\[
	\lim_{n\to \infty} \mu(\vphi_n) 
	= (1-\delta) \mu(\vphi) \geq (1-\delta)M'.
\]
This implies
\[
	\intu{\Es} f_n d\mu \geq 
	\lim_{n\to \infty} \mu(\vphi_n)
	\geq (1-\delta)M'.
\]
By the assumption, 
$\intu{\Es} f_n d\mu \leq \intu{\Es} f d\mu=M$.
Thus, 
\[
	\intu{\Es} f_n d\mu =\lim_{n\to\infty} \intu{\Es} f_n d\mu.
\]
\ \qed 

Next we consider the case with signed coefficient.

\begin{theorem}
	Let $\mu$ be a monotone measure,
	with continuity at $\emptyset$ and from below.
	$\{f_n\}_n$  be an increasing sequence of measurable functions
	converges to $f$.
	Assume that $\Es=\Es^{p\pm},\Es^{p\pm}_\mu$,
	and $\intu{\Es} f_1 d\mu>-\infty$, $\intu{\Es} \cfun_X  d\mu<-\infty$.
	Then
\[
	\lim_{n\to\infty} \intu{\Es} f_n d\mu= \intu{\Es} f d\mu.
\] 
\end{theorem}
Proof.\quad
Set $M=\int f d\mu$.
Let $\eps>0$ be an arbitrary positive number,
and set $M'= M-\eps$ if $M<\infty$.
If $M=\infty$, let $M'$ be any positive number. 
Then, we can select $\vphi \in L(\Es,f)$ with
\[
	 \mu(\vphi) > M'.
\]
By the condition $\int f_1 d\mu>-\infty$,
there exists $\vphi_0\in L(\Es,f_1)$ with $\mu(\vphi_0)>-\infty$.
We give representations for these simple functions as follows.
\[
	\vphi = \sum_k b_k \cfun_{B_k},\quad \vphi_0=\sum_k c_k \cfun_{C_k}.
\]
Define $A^{(\delta)} _n$ for any $\delta>0$ as follows. 
\[
	A^{(\delta)} _n = \{ x| f_n(x)\geq f(x)-\delta \}.
\]
Then, we define simple functions $\{ \vphi_n \}_n$ as follows.
\begin{eqnarray*}
	\vphi_n &=& (\vphi-\delta) \cfun_{ A^{(\delta)}_n} + \vphi_0 \cfun_{ A^{(\delta)\ c}_n} \\
			 &=& \sum_k (b_k-\delta) \cfun_{B_k\cap A^{(\delta)}_n} + \sum_k c_k \cfun_{ C_k\cap A^{(\delta)c}_n} 
\end{eqnarray*}
Then  $\vphi_n \in \Es$. 
Using the definition of $A^{(\delta)}_n$, and the fact $\vphi_0\leq f_1\leq f_n$, 
\[
	\vphi_n (x) =  \vphi \cfun_{ A^{(\delta)}_n}(x) + \vphi_0 \cfun_{ A^{(\delta)c}_n}(x)\leq f_n(x).
\]
Then,
\begin{eqnarray}
	\mu(\vphi_n) 
&=&	 \sum_k (b_k-\delta) \mu(B_k\cap A^{(\delta)}_n) + \sum_k c_k \mu( C_k\cap A^{(\delta)c}_n) \nonumber \\
&=& 	\sum_k b_k \mu(B_k\cap A^{(\delta)}_n)- \delta \sum_k  \mu(B_k\cap A^{(\delta)}_n)\nonumber \\
&&		+\sum_k c_k \mu(B_k\cap A^{(\delta)c}_n). \label{rhs1}
\end{eqnarray}
By Lemma \ref{sum_dct} we have
\[
	\sum_k b_k \mu(B_k\cap A^{(\delta)}_n) \to  \sum_k b_k \mu(B_k), 
\]
\[
	\sum_k c_k \mu(B_k\cap A^{(\delta)c}_n)\to 0,
\]
and
\[
	\sum_k  \mu(B_k\cap A^{(\delta)}_n) \leq \int \cfun_X d\mu.
\] 
Terefore, for large $n$
\[
	\int f_n d\mu > \mu(\vphi_n) \geq \int f d\mu -\delta \int \cfun_X d\mu - 2\eps 
\] 
Then we have $\ds \lim_{n\to\infty} f_n d\mu =\int f d\mu$.
This concludes the proof.\qed

{\it REMARK.} \quad 
Reversing the signatures, the above theorem corresponds a monotone decreasing 
convergence theorem for $\intd{\Es}$.

\subsection{ Monotone decreasing convergence theorems for $\intd{\Es}$}
   
   As we remarked in the previous section,
	when $\Es=\Es^{p\pm},\Es^{p\pm}_\mu$,
	monotone decreasing convergence theorems for $\intd{\Es}$
	are essentially same with monotone increasing convergence theorems
	$\intu{\Es}$. However, 
	when $\Es=\Es^{p+},\Es^{p+}_\mu$, situations are quite different,
	in this section we treat this case.

\begin{lemma} \label{suppToEmpty}
Let $\mu$ be a monotone measure, with continuity at $\emptyset$,
$\Es=\Es^{p+},\Es^{p+}_\mu$,
and $\{A_n\}$ be a decreasing set sequence, with 
\[
	\intd{\Es} f d\mu<\infty, \quad 
	\bigcap_n A_n =\emptyset
\]
Then, we have
\[
	\intd{\Es} f \cfun_{A_n} d\mu \searrow 0.
\]
\end{lemma}
Proof.\quad
$\intd{\Es} f d\mu<\infty$ if and only if there exists $\vphi \in U(\Es, f)$ 
with  $\mu(\vphi)<\infty$.
Thus,
\[
	f\cfun_{A_n}\leq \vphi \cfun_{A_n}
	= \sum_k b_k \cfun_{B_k\cap A_n}
\]
implies 
\[
	\int f\cfun_{A_n} d\mu 
	\leq \mu(\vphi \cfun_{A_n})= \sum_k b_k \mu( B_k\cap A_n).
\]
Using the continuity of $\mu$ at $\emptyset$,
\[
	\mu(B_k\cap A_n) \searrow 0,   \quad (n\to \infty)
\]
for each $k\in\NN$.
By Lemma \ref{sum_dct}, we have, 
\[
	\int f \cfun_{A_n} d\mu \leq \sum_k b_k \mu( B_k\cap A_n) \searrow 0.
\]
\ \qed

\begin{theorem}
Let $\mu$ be a monotone measure, with continuity at $\emptyset$,
$\Es=\Es^{p+},\Es^{p+}_\mu$.
Let $\{ f_n\}$ be a decreasing sequence of measurable functions
converges to $f$.
Assume that $\intd{\Es} f_1 d\mu <\infty$,
we have
\[
	\intd{\Es} f_n d\mu \searrow \int f d\mu .
\]
\end{theorem}
Proof.\quad
By Lemma \ref{decIneq2}. 
\begin{equation}
	\intd{\Es} f_n d\mu \leq \intd{\Es} f_n \cfun_{\{f=0\}} d\mu + \intd{\Es} f_n \cfun_{\{f>0\}}d\mu   \label{sd_div}
\end{equation}
$\intd{\Es} f_1 d\mu<\infty$ implies that there exists $\vphi_0 \in \Es$ with 
$\vphi_0 \geq f_1$ and $\mu(\vphi_0)<\infty$.
Then $\intd{\Es} f_n \cfun_{\{f=0\}} d\mu \to 0$ ($n\to \infty$) as follows.
\begin{eqnarray}
&&			\intd{\Es} f_n \cfun_{\{ f=0\}} d\mu  \nonumber\\
&\leq &	\intd{\Es} f_n \cfun_{\{ f=0\}}\cfun_{\{f_n \leq \eps\vphi_0 \}} d\mu
			+\intd{\Es} f_n \cfun_{\{ f=0\}}\cfun_{\{f_n > \eps\vphi_0 \}} d\mu \nonumber \\ 
&\leq &	\intd{\Es} \eps\vphi_0 \cfun_{\{ f=0\}}\cfun_{\{f_n \leq \eps\vphi_0 \}} d\mu
			+\intd{\Es} f_1 \cfun_{\{ f=0\}} \cfun_{\{f_n > \eps\vphi_0 \}} d\mu. \label{rhs2}
\end{eqnarray}
We remark that $f_1(x)>0$ if $f_n(x)>0$.
Then, for every $x$ with $f(x)=0$,  $f_n(x)\leq \eps\vphi_0(x)$ for large enough $n$.
\[
	\{ x: f_n(x)>\eps\vphi_0(x) \} \searrow \emptyset.
\]
By Lemma \ref{sum_dct}, the second term in (\ref{rhs2}) converges to $\to 0$.
\[
 \mbox{(1 st. term of \ref{rhs2}) } \leq \eps \mu(\vphi_0) \to 0 \quad (\eps \to 0).
\]
Thus $\intd{\Es} f_n \cfun_{\{ f=0\}} d\mu \to 0$ ($n\to \infty$). 

Next we show that $\intd{\Es} f_n \cfun_{\{ f>0\}} d\mu \to \intd{\Es} f  d\mu$ ($n\to \infty$).
Fix any $\delta>0$. Set 
\[
	A^{(\delta)}_n = \{x: f_n(x) \leq (1+\delta)f(x) \}.
\]
Then, using the fact that 
\[
	A^{(\delta)}_n \cap \{x:f(x)>0\} \nearrow \{x:f(x)>0\} \ (n\to \infty),
\]
we have:
\begin{eqnarray*}
&& \intd{\Es} f_n \cfun_{\{f>0\}} d\mu \\
&\leq & \intd{\Es} f_n\cfun_{\{f>0\}} \cfun_{A_n^{(\delta)}} d\mu+\intd{\Es} f_n\cfun_{\{f>0\}} \cfun_{A_n^{(\delta)c}} d\mu \\
&\leq & (1+\delta) \intd{\Es} \ f\ \cfun_{\{f>0\}} \cfun_{A_n^{(\delta)}} d\mu
+\intd{\Es} \ f_1\ \cfun_{\{f>0\} \cap A_n^{(\delta)c}} d\mu \\
&\leq & (1+\delta) \intd{\Es} f d\mu+\intd{\Es} f_1\cfun_{\{f>0\} \cap A_n^{(\delta)c}} d\mu .
\end{eqnarray*}
By Lemma \ref{suppToEmpty}, the second term tends to 0. Therefore,
\[
	\inf_n \intd{\Es} f_n \cfun_{\{ f>0\}} d\mu \leq (1+\delta) \intd{\Es} f d\mu \to   \intd{\Es} f d\mu, \quad (\delta \searrow 0).
\]
The reverse inequality is evident.
\[
	\inf_n \intd{\Es} f_n \cfun_{\{ f>0\}} d\mu =\lim_{n\to \infty} \intd{\Es} f_n \cfun_{\{ f>0\}} d\mu 
	=\intd{\Es} f d\mu
\]
Thus this concludes the proof.\qed 

\subsection{Monotone decreasing convergence theorems for $\intu{\Es }$}
Under some special conditions, $\intu{\Es }$ satisfy the monotone decreasing convergence theorem.
In this section we will give some of them.

First, we consider sub-additive case.
\begin{lemma}\label{panEqCav} 
	Let $\mu$ be a monotone measure with 
\[
	A\cap B=\emptyset \Rightarrow \mu(A \cup B)\leq \mu(A) +\mu(B).
\]
(Such a monotone measure is said to be {\it sub-additive}.)
Set $(\Es_1,\Es_2)$$=(\Es^{p+},\Es^{c+}),$ $(\Es^{p+}_\mu,\Es^{c+}_\mu)$.
Then, for any non-negative measurable function $f$, 
\[
	\intu{\Es_1} f d\mu = \intu{\Es_2} f d\mu. 
\] 
\end{lemma}
Proof.\quad
For any non-negative simple function $\vphi\in \Es_1$,
$\mu(\vphi)$ does not decrease when the corresponding partition is replaced by its refinement.
Moreover, for $\psi \in \Es_2$, we can construct $\vphi'\in\Es_1$ with
$\psi(x)=\vphi'(x)$ as two functions.
The sub-additivity implies also $\mu(\psi)\leq \mu(\vphi')$.
Obviously $\Es_1\subset \Es_2$, and this concludes the proof.\qed

For a sub-additive monotone measure, a Pan integral has the following linearity
(\cite{OuyMes2017}).
This is proved for $\intu{\Es^{p+}}$, however, a similar proof valid for $\intu{\Es^{p+}}_\mu$. 
\begin{theorem}\label{PanLinear}
 (Yao Ouyang, Jun Li, Radko Mesiar \cite{OuyMes2017})\quad
Let $\mu$ be a sub-additive monotone measure.
$f,g$ be non-negative measurable functions,
and $a,b$ be non-negative constants.
Assume that $\Es=\Es^{p+},\Es^{p+}_\mu$.
Then,
\[
	\int (a f + b g ) d\mu = a\int f d\mu + b \int g d\mu.
\]
\end{theorem}
\ \qed

\begin{lemma}\label{prePanDec}
	Let $\mu$ be a monotone measure, 
	$\{A_n\}$ be a decreasing sequence of measurable sets with $A_n \searrow \emptyset$.
	Assume that $\Es=\Es^{p+},\Es^{p+}_\mu$ and
	$\int f d\mu <\infty$. Then, 
\[
	\int f \ \cfun_{A_n} d\mu \searrow 0
\]
\end{lemma}
Proof.\quad
By Lemma \ref{decIneq2} (a), 
\[
	\int f \ d\mu \geq 
	\int f \ \cfun_{A_n^c} d\mu + \int f \ \cfun_{A_n} d\mu.
\]
By Theorem \ref{incCt}, we have $\int f \cfun_{A_n^c} d\mu \to \int f d\mu$.
Then,  consider the limit of $n\to\infty$
\[
	\lim_{n\to\infty} \int f \ \cfun_{A_n} d\mu= 0. 
\]
\ \qed

\begin{theorem}
Let $\mu$ be a sub-additive monotone measure, 
$\{ f_n \}$ be a decreasing sequence of measurable functions.
Assume $\Es = \Es^{p+},\Es^{p+}_\mu,\Es^{c+},\Es^{c+}_\mu$,
and $\intu{\Es} f_1 d\mu <\infty$. Then,
\[
	\intu{\Es} f_n d\mu \searrow \int f d\mu.
\]
\end{theorem}
Proof.\quad
By Lemma \ref{panEqCav}, we prove the theorem for Pan integral.  

Fix any $\delta>0$, and set 
\[
	A_n^{(\delta)} = \{x: f_n(x) \leq f(x)+\delta f_1(x) \}
\]
$f(x)=f_n(x)=0$ when $f_1(x)=0$ since the sequence is non-increasing.
Then, $A_n^{(\delta)} \nearrow X$ ($n\to\infty$) .
Thus,  
\begin{eqnarray*}
\int f_n d\mu 
&=& \int f_n\ \cfun_{A_n^{(\delta)}} + f_n\ \cfun_{A_n^{(\delta)\ c}} d\mu \\
&\leq & \int (f\ +\delta f_1)\ \cfun_{A_n^{(\delta)}} d\mu 
			+ \int f_1\ \cfun_{A_n^{(\delta)\ c}} d\mu  \\
&\leq  & \int f \ \cfun_{A_n^{(\delta)}} d\mu  +\delta \int f_1  d\mu + 
			\int f_1\ \cfun_{A_n^{(\delta)\ c}} d\mu.
\end{eqnarray*}
Then, using Theorem\ref{incCt}, the first term of the above formula
converges to $\int f d\mu$.
Then, using Lemma \ref{prePanDec} ($\delta>0$ is arbitrary small),
\[
	\lim_{n\to \infty} \int f_n d\mu  =\inf_{n} \int f_n d\mu  \leq \int f d\mu.
\]
Using the reverse inequality, which is obvious
\[
	\lim_{n\to \infty} \int f_n d\mu  = \int f d\mu.
\]
\ \qed

Monotone decreasing convergence theorem for Pan integral is valid,
when the limit is  $0$(constant function).

\begin{theorem}
$\Es=\Es^{p+},\Es^{p+}_\mu$,
Let $\mu$ be a monotone measure with the continuity at $\emptyset$.
$\{ f_n \}_n$ be a decreasing sequence of measurable functions.
We assume that $\ds \int \cfun_X d\mu<\infty$,
$\int f_1 d\mu <\infty$, and $ f_n \searrow 0,\quad (n\to \infty)$.
Then, 
\[
	\lim_{n\to\infty} \int f_n d\mu=0.
\]
\end{theorem}
Proof. \quad 
Fix an arbitrary $\delta>0$. Set,
\[
	A_n^{(\delta)}= \{x: f_n(x)>\delta\}
\]
Then,  $A_n^{(\delta)}\searrow \emptyset$. 
Let $\vphi=\sum_k b_k\cfun_{B_k}$ 
be an any simple function in  $L(\Es,f_n)$. 
Then,
\[
	\mu(\vphi) = \sum_k b_k \mu(B_k) =  \sum_{b_k\leq\delta}  b_k \mu(B_k) 
	+ \sum_{b_k> \delta} b_k \mu(B_k). 
\]
When $b_k>\delta$,
using $\ds\inf_{x\in B_k} f(x)\geq b_k>\delta$ and
$B_k\subset A^{(\delta)}_n$,
\begin{eqnarray*}
	\mbox{The right hand side. } 
	&\leq&	\delta \sum_k \mu(B_k) + \int f_n \cfun_{A^{(\delta)}_n} d\mu, \\
	&\leq&	\delta \int \cfun_X d\mu + \int f_1 \cfun_{A^{(\delta)}_n} d\mu.
\end{eqnarray*}
$\mu(A^{(\delta)}_n)\to0$ ($n\to\infty$) since $\mu$ is continuous at  $\emptyset$.
Thus, the second term converges to $0$ by Lemma  \ref{prePanDec},
and 
\[
	\inf_n \int f_n d\mu \leq \delta \int \cfun_X d\mu \to 0,\quad (\delta \searrow 0)
\]
\ \qed


\end{document}